\documentclass[11pt]{article}

\usepackage{amsmath}
\usepackage{amssymb}
\usepackage{amsfonts}
\usepackage{epsfig}

\newcommand{\gap}{\vspace{0.1in}}
\newcommand{\epc}{\hspace{1pc}}
\newcommand{\IE}{{\rm I}\!{\rm E}}
\newcommand{\scub}{{\rm {\scriptstyle UB}}}
\newcommand{\sclb}{{\rm {\scriptstyle LB}}}
\newcommand{\scopt}{{\rm {\scriptstyle opt}}}

\newcommand{\thalf}{{\textstyle \frac{1}{2}}}
\newcommand{\sthalf}{{\scriptstyle \thalf}}
\newcommand{\CVaR}{\mbox{CVaR}}
\newcommand{\VaR}{\mbox{VaR}}
\newcommand{\scCVaR}{\rm {\scriptstyle CVaR}}
\newcommand{\scVaR}{\rm {\scriptstyle VaR}}
\newcommand{\scI}{\rm {\scriptstyle I}}
\newcommand{\scII}{\rm {\scriptstyle II}}
\newcommand{\scIII}{\rm {\scriptstyle III}}

\newtheorem{lemma}{Lemma}[section]
\newtheorem{proposition}[lemma]{Proposition}
\newtheorem{corollary}[lemma]{Corollary}
\newtheorem{theorem}[lemma]{Theorem}

\title{On the Global Minimization of the Value-at-Risk\footnote{Preprint ANL/MCS-P1112-1203}}

\author{Jong-Shi Pang\footnote{Department of Mathematical Sciences,
Rensselaer Polytechnic Institute, Troy, New York 12180-3590,
U.S.A. Email: pang@mts.jhu.edu.  The work of this author's
research was partially supported by the National Science
Foundation under grant CCR-0098013.} \and Sven
Leyffer\footnote{Mathematics and Computer Science Division,
Argonne National Laboratory, Argonne, Illinois 60439.  Supported
by the Mathematical, Information, and Computational Sciences
Division subprogram of the Office of Advanced Scientific
Computing Research, Office of Science, U.S.\ Department of Energy,
under Contract W-31-109-ENG-38.}} 

\begin{document}

\maketitle

\begin{abstract}
In this paper, we consider the nonconvex minimization problem of
the value-at-risk (VaR) that arises from financial risk analysis.
By considering this problem as a special linear program with linear
complementarity constraints (a bilevel linear program to be more
precise), we develop upper and lower bounds for the minimum VaR
and show how the combined bounding procedures can be used to
compute the latter value to global optimality.  A numerical example
is provided to illustrate the methodology.
\end{abstract}

\noindent {\bf Dedication.}  With great pleasure we dedicate
this paper to a respected pioneer of our field, Professor Olvi L.\ Mangasarian,
on the occasion of his 70th birthday.  The two topics of this paper,
LPECs and smoothing methods, are examples of the vast contributions that
Olvi has made in optimization, which have benefited us in many ways 
and which will continue to benefit us
in the future.  Happy 70th birthday, Olvi!

\section{Introduction}

The value-at-risk (VaR) and conditional value-at-risk (CVaR) are
two important risk measures that have been used extensively in recent
years in portfolio selection and in risk analysis.  Whereas the VaR is
closely related to a particular quantile of a random variable, the CVaR
is formally defined and analyzed by Rockafellar and Uryasev in two
papers \cite{RUryasev00,RUryasev02} as a way to alleviate
some of the computational difficulties associated with the optimization
of the VaR.  There is now a substantial literature on the applications
and further developments of these two risk measures; a partial list of
this literature relevant to optimization includes
the papers \cite{AMRUryasev01,KPUryasev02,LMUryasev02,ORuszczynski02,Pflug00,
RUZabarankin02,RUZabarankin03,TUryasev02}. In particular, the paper
\cite{LMUryasev02} presents some CVaR-based algorithms for computing the VaR in
a portfolio selection problem; in spite of their practical efficiency,
however, these algorithms offer no guarantee of global optimality of the 
computed VaR.

Setting aside some criticisms of the VaR mentioned in the literature,
part of which stems from the difficulty associated with the portfolio
selection problem using the VaR criterion, we study the global
optimization problem using a scenario formulation, which is the principal
approach employed in the cited references for solving the (C)VaR minimization
problem.  Specifically, we consider the VaR minimization problem as an
LPEC \cite{Mangasarian98}, a linear program with 
equilibrium constraints, which is a subject pioneered by Mangasarian,
to whom this paper is dedicated.  By exploiting the
special structure of this program, we derive linear
programs whose optimum objective values yield upper and lower bounds for
the optimal VaR.
The bounding procedures are then used in a branch-and-cut
algorithm for computing the latter value to global optimality.
A numerical example is provided to illustrate the algorithm.

An LPEC is a special case of a mathematical program with equilibrium
constraints (MPEC).  Since the publication of the two
monographs \cite{LPRalph96,OKZowe98}, there has been significant
computational advance in numerical methods for solving MPECs;
a partial list of recent references includes
\cite{FJQi99,FLPa98,FPang98,FPang99,FLeyffer02,FLRScholtes02,FTseng02,JRalph00,JRalph03,
Leyffer02,SStohr99,SScholtes00,Scholtes01}.  In spite of such extensive efforts, the
computation of globally optimal solutions to MPECs remains elusive.  While some MPEC
solvers are fairly robust in practice, there is no guarantee that their computed
solutions are globally optimal solutions.  An important reason for this lack of
guarantee for global optimality is the fact that these solvers are all based on local
improvement techniques and no global optimization is incorporated in their
implementation.  As a special MPEC, the VaR minimization problem is amenable
to solution by any one of the (local) methods.  In this paper, we do not stop
with this routine adaptation of the existing MPEC solvers; instead, our goal
is to develop a branch-and-cut algorithm for solving the
minimum VaR problem to global optimality.

\section{The VaR Minimization Problem}

Let $y$ denote an $n$-dimensional random vector whose components represent
the random losses of some financial instruments.  Let $X \subseteq \Re^n$ be
a closed convex set (polyhedral in many practical applications) representing
the set of feasible investments.  For a given $x \in \Re^n$, $z \equiv x^Ty$
is therefore the random loss associated with the investment vector $x$.  For a given
scalar $\beta \in (0,1)$, which denotes a confidence threshold of sustainable loss,
the CVaR and VaR associated with the random variable $z$ is given, as proved
in \cite{RUryasev00,RUryasev02}, by the following two deterministic quantities,
respectively:
\begin{eqnarray*}
\CVaR_{\beta}(x) & \equiv & \min_{m \in \Re} \, \left[ \, m + \frac{1}{1 - \beta} \,
\IE_y ( \, x^Ty - m \, )_+ \, \right], \\ [10pt]
\VaR_{\beta}(x)  & \equiv & \min \{ \, m \, : \, m \, \in \, {\cal M}_{\beta}(x) \, \},
\end{eqnarray*}
where $\IE_y$ denotes the expectation with respect to the random vector $y$, the subscript
plus sign denotes the nonnegative part of a scalar (i.e., the {\sl plus function}
$t_+ \equiv \max(0,t)$), and
${\cal M}_{\beta}(x)$ denotes the set of minimizers in the definition of $\CVaR_{\beta}(x)$.
By the results in the cited references, $\CVaR_{\beta}(x)$ and $\VaR_{\beta}(x)$ are well-defined
finite scalars for very general loss distributions.  Clearly, we have
\[
\CVaR_{\beta}(x) \, = \, \VaR_{\beta}(x) + \frac{1}{1 - \beta} \, \IE_y ( x^Ty - \VaR_{\beta}(x) )_+
\, \geq \, \VaR_{\beta}(x), \epc \forall \, x.
\]

The CVaR and VaR minimization problems are, respectively,
\[
\left\{ \begin{array}{ll}
\mbox{minimize}   & \CVaR_{\beta}(x) \\ [5pt]
\mbox{subject to} & x \, \in \, X
\end{array} \right\} \epc \mbox{and} \epc
\left\{ \begin{array}{ll}
\mbox{minimize}   & \VaR_{\beta}(x) \\ [5pt]
\mbox{subject to} & x \, \in \, X
\end{array} \right\}.
\]
Clearly, the CVaR minimization problem can be cast equivalently as the
following convex program in the joint variable $(m,x)$:
\[ \begin{array}{ll}
\mbox{minimize} & m + \displaystyle{
\frac{1}{1 - \beta}
} \, \IE_y ( \, x^Ty - m \, )_+ \\ [10pt]
\mbox{subject to} & ( \, m,x \, ) \, \in \, \Re \times X.
\end{array} \]
Nevertheless, the VaR problem is {\sl not}
a convex program; this fact is an acknowledged drawback of using the VaR as a criterion in
portfolio selection.  Our main goal in this paper is to develop remedies to this drawback.

\subsection{An LPEC Formulation}

In the rest of the paper, we take $X$ to be a compact polyhedron.  We adopt a scenario
approach to discretize the random vector $y$.  With this approach, the CVaR minimization
problem becomes a linear
program (LP) and the VaR becomes a bilevel linear program, which we reformulate as an
LPEC using the optimality conditions of the lower-level LP.  Specifically,
let $\{y^1, \cdots, y^k\}$ be the finite set of scenario values of $y$, and let
$\{p_1, \cdots, p_k\}$ be the associated probabilities of the respective scenarios, which,
summing to one, are assumed to be all positive.  The discretized CVaR minimization problem is
\[ \begin{array}{ll}
\mbox{minimize} & m + \displaystyle{
\frac{1}{1 - \beta}
} \, \displaystyle{
\sum_{i=1}^k
} \, p_i \, ( \, x^Ty^i - m)_+ \\ [0.2in]
\mbox{subject to} & ( \, m,x \, ) \, \in \, \Re \times X,
\end{array} \]
which is equivalent to the linear program in the variables $(m,x,\tau)$:
\begin{equation} \label{eq:CVaR minimization}
\begin{array}{ll}
\mbox{minimize} & m + \displaystyle{
\frac{1}{1 - \beta}
} \, \displaystyle{
\sum_{i=1}^k
} \, p_i \, \tau_i \\ [0.2in]
\mbox{subject to} & x \, \in \, X \\ [5pt]
\mbox{and} & \left\{ \begin{array}{lll}
\tau_i & \geq & 0 \\ [5pt]
\tau_i & \geq & x^Ty^i - m
\end{array} \right\} \epc \forall \, i \, = \, 1, \ldots, k.
\end{array} \end{equation}

For a given $x \in X$, $\CVaR_{\beta}(x)$ is the minimum objective value of the following
simple LP in the variable $(m,\tau) \in \Re^{1+k}$:
\begin{equation} \label{eq:CVaR LP}
\begin{array}{ll}
\mbox{minimize} & m + \displaystyle{
\frac{1}{1 - \beta}
} \, \displaystyle{
\sum_{i=1}^k
} \, p_i \, \tau_i \\ [0.2in]
\mbox{subject to} & \left\{ \begin{array}{lll}
\tau_i & \geq & 0 \\ [5pt]
\tau_i & \geq & x^Ty^i - m
\end{array} \right\} \epc \forall \, i \, = \, 1, \ldots, k.
\end{array} \end{equation}
By letting $\lambda_i$ denote the dual variable of the $i$th functional constraint
in (\ref{eq:CVaR LP}), the above LP can be solved trivially via its dual:
\begin{equation} \label{eq:dual CVaR LP}
\begin{array}{ll}
\mbox{maximize} & \displaystyle{
\sum_{i=1}^k
} \, \lambda_i \, x^Ty^i \\ [0.2in]
\mbox{subject to} & 0 \, \leq \, \lambda_i \, \leq \, p_i/( 1 - \beta ), \epc \forall
\, i \, = \, 1, \ldots, k \\ [10pt]
\mbox{and} & \displaystyle{
\sum_{i=1}^k
} \, \lambda_i \, = \, 1,
\end{array} \end{equation}
which is a bounded knapsack problem that can in turn be solved by a simple
sorting procedure.  The optimal objective value of either (\ref{eq:CVaR LP}) or
(\ref{eq:dual CVaR LP}) yields CVaR$_{\beta}(x)$; this shows in particular
that CVaR$_{\beta}(x)$ is a convex combination of the portfolio losses
$\{ x^Ty^1, \cdots, x^Ty^k \}$.

In general, by solving
either of the LPs (\ref{eq:CVaR LP}) or (\ref{eq:dual CVaR LP}), we are
not guaranteed to obtain $\VaR_{\beta}(x)$ right away; to obtain the latter value,
we can solve another simple LP in the variable $(m,\tau)$, with $x$ remaining fixed:
\[ \begin{array}{ll}
\mbox{minimize} & m \\ [5pt]
\mbox{subject to} & m + \displaystyle{
\frac{1}{1 - \beta}
} \, \displaystyle{
\sum_{i=1}^k
} \, p_i \, \tau_i \, \leq \, \CVaR(x) \\ [0.2in]
\mbox{and} & \left\{ \begin{array}{lll}
\tau_i & \geq & 0 \\ [5pt]
\tau_i & \geq & x^Ty^i - m
\end{array} \right\} \epc \forall \, i \, = \, 1, \ldots, k,
\end{array} \]
which is simply the problem of finding the least element of the argmin ${\cal M}_{\beta}(x)$.

The optimality conditions of (\ref{eq:CVaR LP}) are
\[ \begin{array}{l}
\left\{ \begin{array}{lllll}
0 & \leq & \tau_i & \perp & \displaystyle{
\frac{p_i}{1 - \beta}
} - \lambda_i \, \geq \, 0 \\ [0.2in]
0 & \leq & \lambda_i & \perp & s_i \, \equiv m + \tau_i - x^Ty^i \, \geq \, 0
\end{array} \right\} \epc \forall \, i \, = \, 1, \ldots, k \\ [0.3in]
\mbox{and} \epc \displaystyle{
\sum_{i=1}^k
} \, \lambda_i \, = \, 1,
\end{array} \]
where the $\perp$ denotes the well-known complementary slackness condition.
Employing these optimality conditions, we can reformulate the VaR minimization
problem as the following linear program with linear complementarity constraints
in the variables $(m,x,\tau,\lambda)$, that is, an LPEC, which in turn is a special
subclass of the class of mathematical programs with equilibrium constraints
\cite{LPRalph96}:
\begin{equation} \label{eq:VaR MPEC}
\begin{array}{ll}
\mbox{minimize}    & m \\ [5pt]
\mbox{subject to } & x \, \in \, X \\ [5pt]
& \left\{ \begin{array}{lllll}
0 & \leq & \tau_i & \perp & \displaystyle{
\frac{p_i}{1 - \beta}
} - \lambda_i \, \geq \, 0 \\ [0.2in]
0 & \leq & \lambda_i & \perp & s_i \, \equiv \, m + \tau_i - x^Ty^i \, \geq \, 0
\end{array} \right\} \epc \forall \, i \, = \, 1, \ldots, k \\ [0.3in]
\mbox{and} & \displaystyle{
\sum_{i=1}^k
} \, \lambda_i \, = \, 1.
\end{array} \end{equation}
As an LPEC, the feasible region of (\ref{eq:VaR MPEC}) is the union
of finitely many polyhedra.  Exploiting its special structure, we
state and prove in the result below that (\ref{eq:VaR MPEC}) attains
a finite minimum objective value.

\begin{proposition} \label{pr:existence VaR} \rm
Let $X$ be a compact polyhedron in $\Re^n$.  The LPEC (\ref{eq:VaR MPEC})
attains a finite minimum objective value.
\end{proposition}

\noindent {\bf Proof.}  Since $X$ is compact by assumption, one can easily
show that $m$ must be bounded below on the feasible region of (\ref{eq:VaR MPEC}).
In fact, if $(m_{\nu},\tau^{\nu},x^{\nu})$ is a sequence of feasible solutions
with $m_{\nu} \to -\infty$, then $\tau^{\nu}_i \to \infty$ for every $i$.  Consequently,
$\lambda^{\nu}_i = p_i/(1-\beta)$; but this contradicts the last constraint, which
requires that the sum of the $\lambda$'s be equal to unity.  \hfill $\Box$

\gap

Let $(m_{\scVaR},x^{\scVaR},\tau^{\scVaR},\lambda^{\scVaR})$ denote an optimal
solution of (\ref{eq:VaR MPEC}).  Note that whereas $m_{\scVaR}$ must be unique,
the triple $(x^{\scVaR},\tau^{\scVaR},\lambda^{\scVaR})$ is not necessarily
so.  Our goal is to compute $m_{\scVaR}$ as best as possible.  Although a
theoretical guarantee of global optimality is not easy to obtain, we derive
valid upper and lower bounds for $m_{\scVaR}$ and develop ways to tighten
these bounds; obviously, when the upper and bounds
coincide, then $m_{\scVaR}$ is obtained.

\section{Upper and Lower Bounds} \label{sec:bounds}

In this section, we develop valid upper and lower bounds for $m_{\scVaR}$.
While upper bounds are not difficult to compute, sharp lower bounds are less
obvious to derive.   We formally describe these bounds in the next two
subsections.  Here, we note that if $m_{\scVaR} \in [ m_{\sclb},m_{\scub} ]$, then
\[
0 \, \leq \, \max \left( \, \frac{m_{\scub} - m_{\scVaR}}{m_{\scVaR}}, \,
\frac{m_{\scVaR} - m_{\sclb}}{m_{\scVaR}} \, \right) \, \leq \, \frac{m_{\scub} - m_{\sclb}}{m_{\sclb}},
\]
which gives relative accuracies of the upper and lower bound values, $m_{\scub}$ and $m_{\sclb}$,
respectively, with respect to the exact minimum VaR $m_{\scVaR}$.

\subsection{LP Upper Bounds} \label{subsec:upper bounding LP}

In essence, upper bounds for $m_{\scVaR}$ are obtained by ``breaking''
the complementary slackness in (\ref{eq:VaR MPEC}) (i.e., restricting
the feasible region) according to a given feasible solution.  Let
$x^0 \in X$ be given.  The scalar
$m_0 \equiv \VaR_{\beta}(x^0)$ provides an upper bound for $m_{\scVaR}$.
(For instance, we may take $x^0 = x^{\scCVaR}$ to be an optimal solution of
the CVaR linear program (\ref{eq:CVaR minimization}).)
We wish to improve on the bound $m_0$ by considering a restriction
of the constraints in (\ref{eq:VaR MPEC}).  Specifically, associated
with the pair $(m_0,x^0)$, let $(\tau^0,\lambda^0)$ satisfy
\[ \begin{array}{l}
\left\{ \begin{array}{lllll}
0 & \leq & \tau_i^0 & \perp & \displaystyle{
\frac{p_i}{1 - \beta}
} - \lambda_i^0 \, \geq \, 0 \\ [0.2in]
0 & \leq & \lambda_i^0 & \perp & s_i^0 \, \equiv m_0 + \tau_i^0 - ( \, x^0 \, )^Ty^i
\, \geq \, 0
\end{array} \right\} \epc \forall \, i \, = \, 1, \ldots, k \\ [0.3in]
\mbox{and} \epc \displaystyle{
\sum_{i=1}^k
} \, \lambda_i^0 \, = \, 1.
\end{array} \]
Define the index sets
\[ \begin{array}{lll}
\alpha_{\tau}^0 & \equiv & \left\{ \, i \, : \, \tau_i^0 \, > \, 0 \, = \,
\displaystyle{
\frac{p_i}{1 - \beta}
} - \lambda_i^0
\right\} \\ [0.2in]
\beta_{\tau}^0 & \equiv & \left\{ \, i \, : \, \tau_i^0 \, = \, 0 \, = \,
\displaystyle{
\frac{p_i}{1 - \beta}
} - \lambda_i^0
\right\} \\ [0.2in]
\gamma_{\tau}^0 & \equiv & \left\{ \, i \, : \, \tau_i^0 \, = \, 0 \, < \,
\displaystyle{
\frac{p_i}{1 - \beta}
} - \lambda_i^0
\right\}
\end{array} \]
and
\[ \begin{array}{lll}
\alpha_{\lambda}^0 & \equiv & \left\{ \, i \, : \, \lambda_i^0 \, > \, 0 \, = \,
m_0 + \tau_i^0 - ( \, x^0 \, )^Ty^i \, \right\} \\ [10pt]
\beta_{\lambda}^0 & \equiv & \left\{ \, i \, : \, \lambda_i^0 \, = \, 0 \, = \,
m_0 + \tau_i^0 - ( \, x^0 \, )^Ty^i \, \right\} \\ [10pt]
\gamma_{\lambda}^0 & \equiv & \left\{ \, i \, : \, \lambda_i^0 \, = \, 0 \, < \,
m_0 + \tau_i^0 - ( \, x^0 \, )^Ty^i \, \right\}.
\end{array} \]
Let $\delta_{\tau}^0$ and $\delta_{\lambda}^0$ be arbitrary subsets of $\beta_{\tau}^0$ and
$\beta_{\lambda}^0$, respectively.  Consider the following linear program in the variables
$(m,x,\tau,\lambda)$:
\begin{equation} \label{eq:restricted VaR MPEC}
\begin{array}{ll}
\mbox{minimize}   & m \\ [5pt]
\mbox{subject to} & x \, \in \, X \\ [5pt]
& \tau_i \, \geq \, 0 \, = \, \displaystyle{
\frac{p_i}{1 - \beta}
} - \lambda_i, \epc \forall \, i \, \in \, \alpha_{\tau}^0 \, \cup \, \delta_{\tau}^0 \\ [0.2in]
& \tau_i \, = \, 0 \, \leq \, \displaystyle{
\frac{p_i}{1 - \beta}
} - \lambda_i, \epc \forall \, i \, \in \, \gamma_{\tau}^0 \, \cup \,
( \, \beta_{\tau}^0 \setminus \delta_{\tau}^0 \, ) \\ [0.2in]
& \lambda_i \, \geq \, 0 \, = \, m + \tau_i - x^Ty^i,
\epc \forall \, i \, \in \, \alpha_{\lambda}^0 \, \cup \, \delta_{\lambda}^0 \\ [10pt]
& \lambda_i \, = \, 0 \, \leq \, m + \tau_i - x^Ty^i,
\epc \forall \, i \, \in \, \gamma_{\lambda}^0 \, \cup \, ( \, \beta_{\lambda}^0 \setminus \delta_{\lambda}^0 \, )
\\ [10pt]
\mbox{and} & \displaystyle{
\sum_{i=1}^k
} \, \lambda_i \, = \, 1,
\end{array} \end{equation}
which is obtained by restricting the complementarity constraints in (\ref{eq:VaR MPEC}) based
on the above index sets.  Obviously, (\ref{eq:restricted VaR MPEC}) is equivalent to a
simplified LP in the variables $(m,\tau)$ only, with the $\lambda$ variable being removed,
that is,
\begin{equation} \label{eq:simplified restricted VaR MPEC}
\begin{array}{ll}
\mbox{minimize}   & m \\ [5pt]
\mbox{subject to} & x \, \in \, X \\ [5pt]
& \tau_i \, \geq \, 0, \epc \forall \, i \, \in \, \alpha_{\tau}^0 \, \cup \, \delta_{\tau}^0 \\ [10pt]
& \tau_i \, = \, 0, \epc \forall \, i \, \in \, \gamma_{\tau}^0 \, \cup \,
( \, \beta_{\tau}^0 \setminus \delta_{\tau}^0 \, ) \\ [10pt]
& m + \tau_i - x^Ty^i \, = \, 0,
\epc \forall \, i \, \in \, \alpha_{\lambda}^0 \, \cup \, \delta_{\lambda}^0 \\ [10pt]
& m + \tau_i - x^Ty^i \, \geq \, 0,
\epc \forall \, i \, \in \, \gamma_{\lambda}^0 \, \cup \, ( \, \beta_{\lambda}^0 \setminus \delta_{\lambda}^0 \, ).
\end{array} \end{equation}
It is clear that $(m_0,x^0,\tau^0)$ is feasible to (\ref{eq:simplified restricted VaR MPEC}).
In general, if $(m,x,\tau)$ is feasible to (\ref{eq:simplified restricted VaR MPEC}),
then $(m,x,\tau,\lambda^0)$ is feasible to (\ref{eq:restricted VaR MPEC}),
and hence to (\ref{eq:VaR MPEC}).  Consequently, (\ref{eq:simplified restricted VaR MPEC})
attains a finite global minimum.  Moreover, if $\left( m_{1+\sthalf},x^1,\tau^1 \right)$
denotes an optimal solution of (\ref{eq:simplified restricted VaR MPEC}), we must have
$m_{1+\sthalf} \in {\cal M}_{\beta}(x^1)$.  Hence, with $m_1 \equiv \VaR_{\beta}(x^1)$,
we have
\[
m_0 \, \geq m_{1+\sthalf} \, \geq \, m_1 \, \geq \, m_{\scVaR}.
\]
One of two cases must occur: (a) $m_0 = m_1$ (no improvement), or (b)
$m_0 > m_1$ (strict improvement).  In case (a), no improvement
is obtained with the particular choice of the pair of index sets $(\delta_{\tau}^0,\delta_{\lambda}^0)$.
One can then try a new pair and solve a new LP (\ref{eq:simplified restricted VaR MPEC}),
hoping to obtain a strictly improved bound for $m_{\scVaR}$.  In case (b),
we can replace $(m_0,x^0)$ by the pair $(m_1,x^1)$ and repeat the
above procedure.  The following result shows that if strict improvement is obtained at each
iteration, then in a finite number of steps the exact minimum VaR is found.

\begin{theorem} \label{th:finite termination} \rm
Let $\{x^{\nu}\} \subset X$ be a sequence of feasible vectors such that for each $\nu$, $x^{\nu+1}$
is obtained from $x^{\nu}$ by solving a certain restricted LP as described above.
If $\VaR_{\beta}(x^{\nu}) > \VaR_{\beta}(x^{\nu+1})$ for every $\nu$, then a finite $\nu_0$
exists such that $\VaR_{\beta}(x^{\nu_0}) = m_{\scVaR}$.
\end{theorem}

\noindent {\bf Proof.}  The feasible region of (\ref{eq:VaR MPEC}) is the union of finitely
many polyhedra, each being the feasible set of (\ref{eq:restricted VaR MPEC}) corresponding
to a particular tuple of index sets $(\alpha_{\tau}^0,\delta_{\tau}^0,\gamma_{\tau}^0,
\alpha_{\lambda}^0,\delta_{\lambda}^0,\gamma_{\lambda}^0)$.  Since
$\VaR_{\beta}(x^{\nu}) > \VaR_{\beta}(x^{\nu+1})$ for every $\nu$, the tuples of
index sets used to produce the sequence $\{x^{\nu}\}$ cannot repeat.  Since there are
only finitely many such tuples of index sets, in generating the sequence $\{x^{\nu}\}$ we must
have encountered all of them; in other words, we must have searched over
the entire feasible region of (\ref{eq:VaR MPEC}).  Consequently, we must have
$\VaR_{\beta}(x^{\nu_0}) = m_{\scVaR}$ for some $\nu_0$.  \hfill $\Box$

\gap

The above result is mainly of theoretical interest because rarely is one so lucky that
strict improvement can be obtained with each trial choice of $(\delta_{\tau}^0,\delta_{\lambda}^0)$.
Notice that the procedure described herein is based on the premise that the set
$\beta_{\tau}^0 \cup \beta_{\lambda}^0$ is nonempty, which means that the pair
$(m_0,x^0)$ is a {\sl degenerate} feasible solution of (\ref{eq:VaR MPEC}),
degenerate with reference to the complementarity conditions.
When $(m_0,x^0)$ is nondegenerate, we will not able to continue the procedure.
Consequently, this is one of the rare instances in mathematical programming where
{\sl degeneracy actually helps}: it enables one to continue the search for
an improvement in a global optimization procedure.

\subsection{NLP Upper Bounds} \label{subsec:upper bounding NLP}

An alternative approach to obtain an upper bound is to form the equivalent nonlinear
program (NLP) of the LPEC (\ref{eq:VaR MPEC})
\begin{equation} \label{eq:VaR NLP}
\begin{array}{ll}
\mbox{minimize}     & m                         \\ [5pt]
\mbox{subject to }  & x \, \in \, X                 \\ [5pt]
                    & \left\{ \begin{array}{llll}
0 & \leq & \tau_i, & \displaystyle{\frac{p_i}{1 - \beta}} - \lambda_i \, \geq \, 0 \\ [0.2in]
0 & \leq & \lambda_i, & s_i \, \equiv \, m + \tau_i - x^Ty^i \, \geq \, 0
\end{array} \right\} \epc \forall \, i \, = \, 1, \ldots, k     \\ [0.3in]
                    & \displaystyle{\sum_{i=1}^k} \, \lambda_i \, = \, 1, \\ [0.2in]
\mbox{and}          & \displaystyle{\sum_{i=1}^k} \,
                    \left[ \, \tau_i \left(\frac{p_i}{1-\beta} -\lambda_i\right)
                    + \lambda_i s_i \, \right] \, \leq \, 0,
\end{array} \end{equation}
and solve this NLP using standard solvers. The last constraint in this problem is
the complementarity constraint. Note that we do not require a lower bound
on the complementarity constraint in (\ref{eq:VaR NLP}) because all terms in this
expression are nonnegative.

It is well known that the NLP (\ref{eq:VaR NLP}) fails the
Mangasarian-Fromovitz constraint qualification (MFCQ) \cite{Mangasarian69,MFromovitz67}
at {\em any\/} feasible point.  This fact implies that the multiplier set of 
(\ref{eq:VaR NLP}) is unbounded, the central path fails to exist, and active 
constraint normals are linearly dependent.  As a consequence, solving MPECs as
NLPs has been commonly regarded as numerically unsafe.
Recently, however, it has been demonstrated that standard NLP solvers can be
employed to solve the equivalent NLPs of MPECs reliably and efficiently.
The convergence of sequential quadratic programming methods to a ``stationary
point'' of an MPEC is analyzed in \cite{AnitM:00b,FLRScholtes02}, and the extension
of interior point methods to MPECs is described in \cite{LiuSun:02,RaghBieg:02a}.
For other related methods, see \cite{FJQi99,FLeyffer02,FLPa98,FPang99,JRalph00,JRalph03},
and the monographs \cite{LPRalph96,OKZowe98}.

Unfortunately, solving the equivalent NLP (\ref{eq:VaR NLP}) does not in itself
guarantee global optimality, despite the practical success of NLP solvers.  The reason
is that the nonconvex nature of the complementarity constraint implies that
NLP solvers may fail to find the global minimum, or even a feasible point.
Nevertheless, NLP solvers have been shown to provide good solutions for many
practical MPECs \cite{FLeyffer02,RaghBieg:02}, and this is the feature we
wish to exploit here.  In fact, for the numerical example reported in
Section~\ref{sec:numerical}, an NLP solver finds a solution, which we show through
additional techniques is a global minimum.  We note that
the latter proof is demonstrated not by NLP but rather by exhibiting an upper
bound for $m_{\scVaR}$ that coincides with a lower bound.

\subsection{LP Lower Bounds}

Upper bounding alone is not enough to verify global optimality of a nonconvex
problem.  In this subsection, we develop some valid lower bounds for $m_{\scVaR}$.
As a first remark, we note that the simple LP relaxation of (\ref{eq:VaR MPEC}) is
\begin{equation} \label{eq:simple LP}
\begin{array}{ll}
\mbox{minimize}    & m \\ [5pt]
\mbox{subject to } & x \, \in \, X \\ [5pt]
& \left\{ \begin{array}{lll}
0 & \leq & \tau_i, \epc \displaystyle{
\frac{p_i}{1 - \beta}
} - \lambda_i \, \geq \, 0 \\ [0.2in]
0 & \leq & \lambda_i, \epc m + \tau_i - x^Ty^i \, \geq \, 0
\end{array} \right\} \epc \forall \, i \, = \, 1, \ldots, k \\ [0.3in]
\mbox{and} & \displaystyle{
\sum_{i=1}^k
} \, \lambda_i \, = \, 1,
\end{array} \end{equation}
which does not have a finite optimal solution because we can make $m$ tend
to $-\infty$ with each $\tau_i \to \infty$.  Therefore, we need to tighten
this relaxation.  The following lemma gives a preliminary lower bound
for $m_{\scVaR}$.

\begin{lemma} \label{lm:lower bound for VaR} \rm
For any feasible tuple $(m,x,\tau,\lambda)$ to (\ref{eq:VaR MPEC}), an index
$i$ exists such that $m \geq x^Ty^i$ for at least one index $i$.  Consequently,
\[
m \, \geq \, \min_{1 \leq j \leq k} \, \min_{x \in X} \, x^Ty^j \, \equiv \, \underline{m}.
\]
\end{lemma}

\noindent {\bf Proof.}  Let $(m,x,\tau,\lambda)$ be an arbitrary feasible tuple
to (\ref{eq:VaR MPEC}).  We must have
\[
\tau_i \, = \, \max(0, x^Ty^i - m), \epc \forall \, i \, = \, 1, \ldots, k.
\]
From the first complementarity constraint in (\ref{eq:VaR MPEC}), we obtain
\begin{equation} \label{eq:first complementarity}
\tau_i \, \lambda_i \, = \, \frac{p_i}{1 - \beta} \, \tau_i,
\end{equation}
which, when used in the second complementarity constraint, yields
\begin{equation} \label{eq:second complementarity}
0 \, = \, m \, \lambda_i + \frac{p_i}{1 - \beta} \, \tau_i - \lambda_i \, x^Ty^i.
\end{equation}
Since the sum of the $\lambda_i$ is equal to unity, we deduce
\begin{equation} \label{eq:formula for m}
m \, = \, \sum_{j=1}^k \, \left[ \, \lambda_j \, x^Ty^j - \frac{p_j}{1 - \beta} \, \tau_j \, \right].
\end{equation}
which yields
\[
m + \sum_{j=1}^k \, \frac{p_j}{1 - \beta} \, \max(0, x^Ty^j - m) \, = \,
\sum_{j=1}^k \, \lambda_j \, x^Ty^j.
\]
If no index $i$ exists such that $x^Ty^i \leq m$, then $m < x^Ty^j$ for all $j$,
and the above identity yields
\[
m \, = \, \frac{1 - \beta}{\beta} \, \displaystyle{
\sum_{j=1}^k
} \, \left( \, \displaystyle{
\frac{p_j}{1 - \beta}
} - \lambda_j \, \right) \, x^Ty^j,
\]
which shows that $m$ is a convex combination of the family $\{ x^Ty^1, \cdots, x^Ty^k \}$.
This is a contradiction.  The last assertion of the lemma is obvious.  \hfill $\Box$.

\gap

In essence, the lower bounding procedure described below aims at removing the
three nonlinear terms $m\lambda_i$, $\tau \lambda_i$, and $\lambda_i x$ in
(\ref{eq:first complementarity}) and (\ref{eq:second complementarity}),
which are the result of the complementarity constraints, while maintaining
some form of these two equations.  It turns out that the first two nonlinear terms
can be completely removed through some suitable substitution, whereas the third one cannot.
The relaxation of (\ref{eq:VaR MPEC}) then employs a single variable $z^i$ to
substitute for $\lambda_i x$ and to remove the identity $z^i = \lambda_i x$
when $\lambda_i$ is strictly between its lower and upper bounds.  Note that the change
of variables implies
\[
x \, = \, \sum_{i=1}^k \, z^i.
\]
Furthermore, if $x \geq 0$ in the set $X$, it follows that
$0 \leq z^i \leq (p_i/(1-\beta)) x$ for all $i$.  More generally, if $|x| \leq a$ for
all $x \in X$, where $a$ is a given nonnegative vector and the absolute sign is meant
componentwise, then $|z^i| \leq (p_i/(1-\beta)) a$ for all $i$.

From the identity (\ref{eq:second complementarity}), we deduce that for any feasible
tuple $(m,x,\tau,\lambda)$ to (\ref{eq:VaR MPEC}),
\[ \begin{array}{lll}
m \, \geq \, 0 & \Rightarrow & \forall \, i, \, \left[ \,
0 \, \leq \,  \, ( \, z^i \, )^Ty^i -
\displaystyle{
\frac{p_i}{1 - \beta}
} \, \tau_i \, \leq
\, \displaystyle{
\frac{p_i}{1 - \beta}
} \, m \, \right] \\ [0.3in]
m \, \leq \, 0 & \Rightarrow & \forall \, i, \, \left[ \,
0 \, \geq \,  \, ( \, z^i \, )^Ty^i -
\displaystyle{
\frac{p_i}{1 - \beta}
} \, \tau_i \, \geq
\, \displaystyle{
\frac{p_i}{1 - \beta}
} \, m \, \right] .
\end{array} \]
Assume for the moment that $m \geq 0$.
This gives rise to the following LP relaxation of (\ref{eq:VaR MPEC})
(we assume that $|x| \leq a$ for all $x \in X$):
\begin{equation} \label{eq:LP relax}
\begin{array}{ll}
   \mbox{minimize}  & m                                 \\ [5pt]
   \mbox{subject to}    & x \, \equiv \, \displaystyle{\sum_{j=1}^k} \, z^j \, \in \, X,\\ [0.25in]
            & m \, = \, \displaystyle{\sum_{j=1}^k} \, \left[ \, ( \, z^j \, )^Ty^j
                - \frac{p_j}{1 - \beta} \, \tau_j \, \right]        \\ [0.3in]
            & \left\{ \begin{array}{r}
    s_i \, \equiv \, m + \tau_i - x^Ty^i \, \geq 0              \\ [10pt]
    0 \, \leq \, ( \, z^i \, )^Ty^i - \displaystyle{\frac{p_i}{1 - \beta}} \,
        \tau_i \, \leq \, \frac{p_i}{1 - \beta} \, m            \\ [0.3in]
    \tau_i \, \geq \, 0                             \\ [10pt]
    | \, z^i \, | \, \leq \, \displaystyle{\frac{p_i}{1 - \beta}} \, a
            \end{array} \right\}, \epc \forall \, i \, = \, 1, \ldots, k.
\end{array}
\end{equation}
We are also interested in investigating the behavior of this lower bound if we branch
on a disjunction.  It follows from (\ref{eq:VaR MPEC}) that there are three possible
branches for $\lambda_i$. Each branch in turn gives rise to a particular implication:
\[ \begin{array}{lllllllll}
0 & = & \lambda_i & & & \Rightarrow & [ \, z^i \, = \, 0 & \mbox{and} & \tau_i \, = \, 0 \, ]
\\ [10pt]
& & \lambda_i & = & \displaystyle{\frac{p_i}{1 - \beta}} & \Rightarrow &
\left[
z^i \, = \, \displaystyle{\frac{p_i}{1 - \beta}} \, x \, \right. & \mbox{and} &
\left. s_i \, = 0 \, \right] \\ [0.3in]
0 & < & \lambda_i & < &  \displaystyle{\frac{p_i}{1 - \beta}} & \Rightarrow &
[ \, \tau_i \, = \, 0 & \mbox{and} & s_i \, = \, 0 \, ].
\end{array} \]
Based on the above implications, we define three LPs by adding the
implications to the lower bounding LP (\ref{eq:LP relax}), respectively. In
what follows, $i_0$ is a fixed but arbitrary index in $\{1, \ldots, k\}$.
\begin{description}
 \item[(LP$_{\scI,i_0}^+$)] This corresponds to the case where $\lambda_{i_0} = 0$ and $m \geq 0$
and consists of the LP relaxation (\ref{eq:LP relax}) with the following additional constraints:
\begin{equation} \label{eq:LP relaxed I}
    z^{i_0} \, = \, 0 \epc \mbox{and} \epc \tau_{i_0} \, = \, 0.
\end{equation}
\item[(LP$_{\scII,i_0}^+$)] This corresponds to the case where $\lambda_{i_0} = p_{i_0}/(1-\beta)$
and $m \geq 0$ and consists of the LP relaxation (\ref{eq:LP relax}) with the following
additional constraints:
\begin{equation} \label{eq:LP relaxed II}
 z^{i_0} \, = \, \frac{p_{i_0}}{1 - \beta} \, x \, \epc \mbox{and} \epc
s_{i_0} \, = \, 0.
\end{equation}
\item[(LP$_{\scIII,i_0}^+$)] This corresponds to the case where
$\lambda_{i_0} \in (0,p_{i_0}/(1-\beta))$ and $m \geq 0$ and consists of the LP relaxation
(\ref{eq:LP relax}) with the following additional constraints:
\begin{equation} \label{eq:LP relaxed III}
 \tau_{i_0} \, = \, s_{i_0} \, =  \, 0.
\end{equation}
\end{description}
Let LP$_{\scI,i_0}^{+,\scopt}$, LP$_{\scII,i_0}^{+,\scopt}$, and LP$_{\scIII,i_0}^{+,\scopt}$
denote the optimal objective values of the above three LPs, respectively.
Consistent with a standard convention in optimization,
we define the minimum objective value of an infeasible LP to be $\infty$.
The next result summarizes the fundamental role of the above LPs for
solving the VaR minimization problem (\ref{eq:VaR MPEC}).

\begin{proposition} \label{pr:VaR LPEC relaxed}
\rm
For any index $i_0$, the following five statements (a)--(e) are valid.
\begin{description}
\item[\rm (a)] If $(m,x,\tau,\lambda)$ is feasible to (\ref{eq:VaR MPEC}) and $m \geq 0$, then
with $z^i \equiv \lambda_i x$ for all $i$, the tuple
$(m,x,z,\tau)$ is feasible to (\ref{eq:LP relax}). It also satisfies the additional
constraints  (\ref{eq:LP relaxed I}) if $\lambda_{i_0} = 0$,
(\ref{eq:LP relaxed II}) if $\lambda_{i_0} = p_{i_0}/(1-\beta)$, and
(\ref{eq:LP relaxed III}) if $0 < \lambda_{i_0} < p_{i_0}/(1-\beta)$.
\item[\rm (b)] If $m_{\scVaR} \geq 0$, then at least one of the three LPs obtained
by adding to (\ref{eq:LP relax}) the cuts (\ref{eq:LP relaxed I}), or
(\ref{eq:LP relaxed II}), or (\ref{eq:LP relaxed III}), must be feasible and, hence,
solvable; in this case,
\begin{equation} \label{eq:VaR lower bound}
m_{\scVaR} \, \geq \,
\min \left( \mbox{ LP}_{\scI,i_0}^{+,\scopt}, \mbox{ LP}_{\scII,i_0}^{+,\scopt},
\mbox{ LP}_{\scIII,i_0}^{+,\scopt} \right) .
\end{equation}
\item[\rm (c)] If $m_{\scub} \geq m_{\scVaR} \geq 0$ and $\infty > $ LP$_{\scI,i_0}^{+,\scopt} > m_{\scub}$,
then for any optimal solution $(x^{\scVaR},\tau^{\scVaR},\lambda^{\scVaR})$ of
(\ref{eq:VaR MPEC}), we must have $\lambda_{i_0}^{\scopt} > 0$, and thus,
$s_{i_0}^{\scVaR} \, \equiv \tau_{i_0}^{\scVaR} + m_{\scVaR} - (x^{\scVaR})^Ty^{i_0} = 0$.
\item[\rm (d)] If $m_{\scub} \geq m_{\scVaR} \geq 0$ and $\infty > $ LP$_{\scII,i_0}^{+,\scopt} > m_{\scub}$,
then for any optimal solution $(x^{\scVaR},\tau^{\scVaR},\lambda^{\scVaR})$ of
(\ref{eq:VaR MPEC}), we must have $\lambda_{i_0}^{\scVaR} < p_{i_0}/(1 - \beta)$, and thus,
$\tau_{i_0}^{\scVaR} = 0$.
\item[\rm (e)] If $m_{\scub} \geq m_{\scVaR} \geq 0$ and $\infty > $ LP$_{\scIII,i_0}^{+,\scopt} > m_{\scub}$,
then for any optimal solution $(x^{\scVaR},\tau^{\scVaR},\lambda^{\scVaR})$ of
(\ref{eq:VaR MPEC}), we must have $\lambda_{i_0}^{\scVaR} = 0$ or $\lambda_{i_0}^{\scVaR} = p_{i_0}/(1 - \beta)$.
\end{description}
\end{proposition}

\noindent {\bf Proof.}  Part (a) does not require a proof.  For part (b), we need only to
prove the bound (\ref{eq:VaR lower bound}).  Let $(x^{\scVaR},\tau^{\scVaR},\lambda^{\scVaR})$
be an arbitrary optimal solution of (\ref{eq:VaR MPEC}) corresponding to $m_{\scVaR}$.
the tuple $(m_{\scVaR},x^{\scVaR},z^{\scVaR},\tau^{\scVaR})$, where
$z^{\scVaR,i} \equiv \lambda_i^{\scVaR} x^{\scVaR}$, is feasible to the one of the
three LPs formed from (\ref{eq:LP relax}) plus (\ref{eq:LP relaxed I}), or
(\ref{eq:LP relaxed II}), or (\ref{eq:LP relaxed III}).
Hence (\ref{eq:VaR lower bound}) follows readily.  To prove (c), one need only note that
if $\lambda_{i_0}^{\scVaR} = 0$, then $(m_{\scVaR},x^{\scVaR},z^{\scVaR},\tau^{\scVaR})$, where
$z^{\scVaR,i} \equiv \lambda_i^{\scVaR} x^{\scVaR}$, is feasible to (\ref{eq:LP relaxed I});
hence $m_{\scVaR} \geq \mbox{LP}_{\scI,i_0}^{+,\scopt}$, which easily yields a contradiction.
The proof of (d) and (e) is similar and not repeated.  \hfill $\Box$

\gap

We can similarly set up three other LPs to handle the case $m_{\scVaR} \leq 0$.  It suffices
to reverse the inequality signs in
\[
0 \, \leq \, ( \, z^i \, )^Ty^i - \displaystyle{
\frac{p_i}{1 - \beta}
} \, \tau_i \, \leq \, \frac{p_i}{1 - \beta} \, m
\]
and use instead
\begin{equation} \label{eq:reversed inequalities}
0 \, \geq \, ( \, z^i \, )^Ty^i - \displaystyle{
\frac{p_i}{1 - \beta}
} \, \tau_i \, \geq \, \frac{p_i}{1 - \beta} \, m.
\end{equation}
Letting LP$_{\scI,i_0}^{-,\scopt}$, LP$_{\scII,i_0}^{-,\scopt}$, and LP$_{\scIII,i_0}^{-,\scopt}$
denote the optimal objective values of the resulting LPs, respectively, we can obtain a result
similar to Proposition~\ref{pr:VaR LPEC relaxed}.  Combining these two results, we arrive at a
desired lower bound for $m_{\scVaR}$.

\begin{corollary} \label{co:VaR lower bound}
\rm It holds that
\[
m_{\scVaR} \geq \min \left\{
\max_{1 \leq j \leq k} \min \left( \mbox{ LP}_{\scI,j}^{+,\scopt}, \mbox{ LP}_{\scII,j}^{+,\scopt},
\mbox{ LP}_{\scIII,j}^{+,\scopt} \right), \
\max_{1 \leq j \leq k} \min \left( \mbox{ LP}_{\scI,j}^{-,\scopt}, \mbox{ LP}_{\scII,j}^{+,\scopt},
\mbox{ LP}_{\scIII,j}^{+,\scopt} \right)  \right\}.
\]
\hfill $\Box$
\end{corollary}

The practical value of the cuts (\ref{eq:LP relaxed I}), (\ref{eq:LP relaxed II}), and
(\ref{eq:LP relaxed III}), and their analogs with the reverse inequalities
(\ref{eq:reversed inequalities}) built in, lies in their ability to improve the
lower bound obtained from (\ref{eq:LP relax}) making it easier to fathom nodes
in the branch-and-cut framework that is described in Section~\ref{sec:bb}.

\subsection{Convex Hull Relaxations}

Alternative lower bounds can be derived by observing that
the only difference between the simple LP (\ref{eq:simple LP})
and the LPEC (\ref{eq:VaR MPEC}) is the absence of the
complementarity constraint. Thus, to tighten the
former LP relaxation, we form a linear relaxation of the
complementarity constraint,
\begin{eqnarray}
   0 & = & \sum_{i=1}^k \left\{ \tau_i \frac{p_i}{1-\beta} - \tau_i \lambda_i + \lambda_i \tau_i
                   + \lambda_i m - \lambda_i x^T y^i \right\}   \nonumber   \\ [5pt]
     & = & m + \sum_{i=1}^k \left\{ \tau_i \frac{p_i}{1-\beta} - \lambda_i x^T y^i \right\} ,
                                    \label{eq:compl}
\end{eqnarray}
where we have used the fact that $\sum \lambda_i = 1$. Observe that the only
nonlinear term in this expression is given by $\lambda_i x^T y^i$.  Next we
show how to construct the convex hull relaxation of this constraint.

We introduce new linear variables $\gamma_i = x^T y^i$ and then replace the nonlinear
terms $\lambda_i \gamma_i$ by $w_i$. Since $w_i = \lambda_i \gamma_i$ is a simple
bilinear expression, we can strengthen this LP relaxation by adding the convex
hull of $w_i = \lambda_i \gamma_i$. Let $L_i$ and $U_i$ be valid lower and upper
bounds on $\gamma_i$, respectively (these can be obtained by solving $2k$
LPs for instance); the convex hull of $w_i = \gamma_i \lambda_i$ is then given by
\begin{equation}
\left. \begin{array}{rcl}
w_i & \geq & L_i \, \lambda_i                          \\ [10pt]
w_i & \geq & \displaystyle{
\frac{p_i}{1-\beta}
} \, \gamma_i + U_i \, \lambda_i - \displaystyle{
\frac{p_i}{1-\beta}
} \, U_i \\ [0.2in]
w_i & \leq & U_i \, \lambda_i \\ [10pt]
w_i & \leq & \displaystyle{
\frac{p_i}{1-\beta}
} \, \gamma_i + L_i \, \lambda_i - \displaystyle{
\frac{p_i}{1-\beta}
} L_i ,
\end{array} \right\} ;
\end{equation}
see, for example \cite{TawaSahi:02}. This gives rise to the LP relaxation
\begin{equation} \label{eq:convex relax}
\begin{array}{ll}
\mbox{minimize}      & m                             \\ [5pt]
\mbox{subject to }   & x \, = \displaystyle{
\sum_{i=1}^k
} \, z^i \, \in \, X  \\ [0.2in]
& \left\{ \begin{array}{lll}
0 & \leq & \tau_i, \epc \displaystyle{\frac{p_i}{1 - \beta}} - \lambda_i \, \geq \, 0 \\ [0.2in]
0 & \leq & \lambda_i, \epc m + \tau_i - x^Ty^i \, \geq \, 0
\end{array} \right\} \epc \forall \, i \, = \, 1, \ldots, k         \\ [0.3in]
& \displaystyle{\sum_{i=1}^k} \, \lambda_i \, = \, 1   \\ [0.2in]
& \gamma_i \, = \, x^T y^i \, \in \, [ \, L_i,U_i \, ], \epc \forall \, i \, = \, 1, \ldots, k
\\ [10pt]
& \left. \begin{array}{rcl}
w_i & \geq & L_i \, \lambda_i \\ [10pt]
w_i & \geq & \displaystyle{
\frac{p_i}{1-\beta}
} \, \gamma_i + U_i \, \lambda_i - \displaystyle{
\frac{p_i}{1-\beta}
} \, U_i \\ [0.2in]
w_i & \leq & U_i \, \lambda_i \\ [10pt]
w_i & \leq & \displaystyle{
\frac{p_i}{1-\beta}
} \, \gamma_i + L_i \, \lambda_i - \displaystyle{
\frac{p_i}{1-\beta}
} L_i ,
\end{array} \right\} \epc \forall \, i \, = \, 1, \ldots, k \\ [0.7in]
\mbox{and} & \displaystyle
0 \, = \, m + \displaystyle{
\sum_{i=1}^k
} \, \left\{ \, \tau_i \, \displaystyle{
\frac{p_i}{1-\beta}
} - w_i \, \right\}.
\end{array}
\end{equation}
Since this LP includes the convex hull relaxation of $w_i = \gamma_i \lambda_i$, it follows
that the LP is bounded whenever the original LPEC (\ref{eq:VaR MPEC}) is bounded.

The two LP bounds (\ref{eq:LP relax}) and (\ref{eq:convex relax}) are nondominating;
the quality of the bounds differs from problem instance to problem instance.
This is confirmed by the numerical example in Section~\ref{sec:numerical}, where,
in one case, (\ref{eq:LP relax}) yields a sharper lower bound, and in the other case,
it is (\ref{eq:convex relax}) that yields a better bound.
Next, we show how the bounds (\ref{eq:LP relax}) and (\ref{eq:convex relax})
can be employed within a branch-and-cut framework to prove the optimality of
a candidate solution of the LPEC (\ref{eq:VaR MPEC}).

\section{Verifying Optimality by Branch-and-Cut} \label{sec:bb}

We briefly outline our approach to prove the optimality of a given candidate
solution to the MPEC (\ref{eq:VaR MPEC}). There are two main ideas: the first is
to construct as small as possible a branch-and-bound tree corresponding to a given
candidate solution, and the second is to exploit the logical implications from
the complementarity constraint to strengthen the LP relaxation as in
Proposition~\ref{pr:VaR LPEC relaxed}.
Section~\ref{sec:numerical} shows how the approach works for a numerical example.

The use of branch-and-bound to solve MPECs is not new. It has been used
in \cite{BarJF:88} to solve some bilevel convex programs.  However, the
scheme proposed here employs the special bounds derived in the preceding section
that are tailored to the minimum VaR problem and that are used to define cuts
that restrict the feasible region of the problem.

In general, the LPEC is initially solved with the
complementarity constraint relaxed.  If this problem yields a solution that is
complementary, then it is also optimal. Otherwise there exists a complementarity
that is violated and we can branch on this complementarity. Branching introduces
two (or three in our case) child problems where the complementarity is broken.
The procedure continues to solve relaxations and branch until an LPEC feasible
solution is found, a problem is infeasible, or its solution is dominated by an
upper bound. This process is best envisioned as a tree search where nodes
correspond to LP relaxations and edges correspond to branches.

Unfortunately, searching the entire branch-and-bound tree is likely to be inefficient.
Instead, we will exploit the branch-and-cut methodology to construct the smallest
tree that can be used to establish optimality of a given feasible solution.
Specifically, let $(m_*,x^*,\lambda^*,\tau^*)$ be a feasible
point of the LPEC (\ref{eq:VaR MPEC}) whose optimality we wish to prove. Note
that $m_*$ is an upper bound on the optimal value of the LPEC. We
then perform several rounds of bound tightening to fix complementary expressions
by solving LP relaxations.  For instance, if we postulate that
$\lambda_i^* = p_i/(1-\beta)$, then we solve the two remaining
LP relaxations in the disjunction. If they produce bounds that are larger than
the given upper bound $m_*$, then we can fix $\lambda_i^* = p_i/(1-\beta)$.
Similar conclusions are possible for the other bounds and variables.
This generates one hopes a short branch-and-bound tree.  The numerical
example presented next illustrates the idea.

\section{Numerical Example} \label{sec:numerical}

The numerical example has the following data:
$n = 3$, $k = 27$, $\beta=0.9$, $p_i = 1/27$ for all $i$,
\[
X \, \equiv \, \left\{ \, x \, \in \, \Re_+^n \, : \,
\displaystyle{
\sum_{i=1}^n
} \, x_i \, = \, 1, \
\displaystyle{
\sum_{i=1}^n
} \, r_i \, x_i \, \geq \, f \, \right\},
\]
where $r \equiv (-1/3,2/3,-1)$ and $f = 1/10$.  The vectors $y^i$ are generated
as follows.  We generate three vectors
\[
\left[ \begin{array}{c}
d^1 \\ [5pt]
d^2 \\ [5pt]
d^3
\end{array} \right] \, \equiv \, \left[ \begin{array}{rrr}
5 & 0 & -6 \\ [5pt]
7 & 0 & -5 \\ [5pt]
2 & 0 & -5
\end{array} \right];
\]
and then we set
\[
y^j_1 \, = \, d^{\, 1}(i1), \epc y^j_2 \, = \, d^{\, 2}(i2), \epc y^j_3 \, = \, d^{\, 3}(i3),
\]
where $j = 6(i1-1)+3(i2-1)+i3$ for $i1,i2,i3=1,2,3$.

\subsection{Finding an Upper Bound}

We first solve the CVaR LP
(\ref{eq:CVaR LP}) and obtain the solution: CVaR = $5.0644$,
$m_{\scCVaR} = 4.8613$, $x^{\scCVaR} = (0.1097,0.6161,0.2742)$, $\tau_1 > 0$,
$\tau_i = 0$ for all $i \geq 2$, $s_i = 0$ for $i = 1, 2, 10$, and
$s_i > 0$ for all other $i$.  The value $m_{\scCVaR}$ is then verified to be
the least element of ${\cal M}_{\beta}(x^{\scCVaR})$.  Based on the
pair $(m_{\scCVaR},x^{\scCVaR})$, we solve the LP (\ref{eq:restricted VaR MPEC})
by setting
\[ \begin{array}{ll}
\alpha_{\tau}^0 \, \cup \, \delta_{\tau}^0 \, = \, \{ 1,2 \}, &
\gamma_{\tau}^0 \, \cup \, (\beta_{\tau}^0 \setminus \delta_{\tau}^0) \, = \,
\{1, \ldots, 27\} \setminus \alpha_{\tau}^0, \\ [10pt]
\alpha_{\lambda}^0 \, \cup \, \delta_{\lambda}^0 \, = \, \{ 1,2,10 \}, &
\gamma_{\lambda}^0 \, \cup \, (\beta_{\lambda}^0 \setminus \delta_{\lambda}^0)
\, = \, \{1, \ldots, 27\} \setminus \alpha_{\lambda}^0.
\end{array} \]
The optimal solution for this LP is as follows: $m_1 = 4.2652$; $\tau_1$ and $\tau_2$
are both positive, and the remaining $\tau_i$ are zero; $s_i = 0$ for
$i = 1, 2, 3, 10$, and the other $s_i$ are all positive.  This solution yields
$\lambda_1 = \lambda_2$ at their common upper bound, which is 0.37037,
implying that $\lambda_3 \geq 0$ and $\lambda_{10} \geq 0$ satisfy
$\lambda_3 + \lambda_{10} = 1 - 2*0.37037 = .26926$.

The pair $(m_{\scCVaR},x^{\scCVaR})$ also belongs to several other pieces
of the feasible region of (\ref{eq:VaR MPEC}). For example, one such piece
corresponds to the index set partitions:
\[ \begin{array}{ll}
\alpha_{\tau}^0 \, \cup \, \delta_{\tau}^0 \, = \, \{ 1,10 \}, &
\gamma_{\tau}^0 \, \cup \, (\beta_{\tau}^0 \setminus \delta_{\tau}^0) \, = \,
\{1, \ldots, 27\} \setminus \alpha_{\tau}^0, \\ [10pt]
\alpha_{\lambda}^0 \, \cup \, \delta_{\lambda}^0 \, = \, \{ 1,2,10 \}, &
\gamma_{\lambda}^0 \, \cup \, (\beta_{\lambda}^0 \setminus \delta_{\lambda}^0)
\, = \, \{1, \ldots, 27\} \setminus \alpha_{\lambda}^0;
\end{array} \]
nevertheless, solving the associated LPs (\ref{eq:restricted VaR MPEC}) on these
other pieces does not yield a lower objective value than $4.2652$.
The same upper bound $m_{\scub} = 4.2652$ is also obtained by solving the
single equivalent NLP (\ref{eq:VaR NLP}).

\subsection{Branching to Verify Global Optimality} \label{subsec:branching optimality}

Our next task is to determine whether the value $m_{\scub} = 4.2652$ is globally optimal.
First, we verify that $m_{\scVaR} \geq 0$ by solving the LP
\[ \begin{array}{ll}
\mbox{minimize} & m \\ [5pt]
\mbox{subject to} & x \, \equiv \, \displaystyle{
\sum_{j=1}^k
} \, z^j \, \in \, X, \\ [0.25in]
& m \, = \, \displaystyle{
\sum_{j=1}^k
} \, \left[ \, ( \, z^j \, )^Ty^j - \frac{p_j}{1 - \beta} \, \tau_j \, \right]  \\ [0.3in]
& \left\{ \begin{array}{r}
s_i \, \equiv \, m + \tau_i - x^Ty^i \, \geq 0 \\ [10pt]
0 \, \geq \, ( \, z^i \, )^Ty^i - \displaystyle{
\frac{p_i}{1 - \beta}
} \, \tau_i \, \geq \, \frac{p_i}{1 - \beta} \, m \\ [0.3in]
\tau_i \, \geq \, 0, \epc 0 \, \leq \, z^i \, \leq \,
\displaystyle{
\frac{p_i}{1 - \beta}
} \, x
\end{array} \right\}, \epc \forall \, i \, = \, 1, \ldots, k.
\end{array} \]
This LP is infeasible, and we can therefore conclude that $m_{\scVaR} \geq 0$.

    \begin{figure}[htb]
    \begin{center}
       \epsfig{ file = 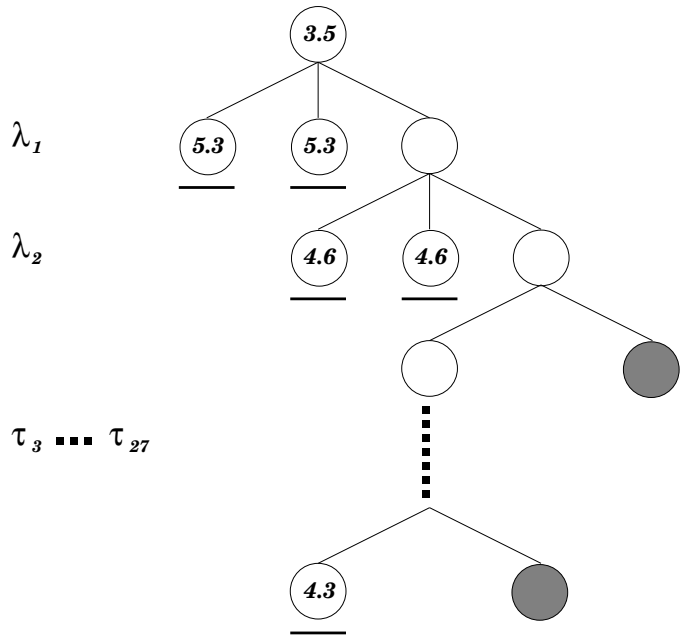 , width = 3.5in }
       \caption{Branch-and-bound tree for the numerical example\label{F-tree}}
    \end{center}
    \end{figure}

Figure~\ref{F-tree} shows the branch-and-bound tree that we construct for
this example.  Each node corresponds to an LP relaxation, with additional
constraints included according to (\ref{eq:LP relaxed I}), (\ref{eq:LP relaxed II}),
or (\ref{eq:LP relaxed III}).  The root node shows the value of the LP relaxation
(\ref{eq:LP relax}) (all values are rounded to two digits).  The variable names
on the left indicate the branching variable. For each $\lambda_i$ there are
three branches (which are ordered from left to write as
$\lambda_i=0, \; \lambda_i \in (0,p_i/(1-\beta)), \; \lambda_i = ,p_i/(1-\beta)$).
For each $\tau_i$ there are two branches (from left to right $\tau_i=0$, and
$\tau_i > 0$).

The lower bounds alone do not allow us to conclude optimality of the
candidate solution.  Hence, we start the construction of the
branch-and-bound tree by proving that $\lambda_1$ and $\lambda_2$
must be at their upper bounds at an optimal solution.  For this purpose,
we solve two LPs by adding the cuts (\ref{eq:LP relaxed I})
and (\ref{eq:LP relaxed III}) to (\ref{eq:LP relax}). These LPs have an
optimal value of $5.3$, which is larger than $m_{\scub}$; we can therefore
consider those nodes as fathomed.  This is illustrated in the tree in
Figure~\ref{F-tree} by the
bold horizontal lines under the node. Hence, we can fix $\lambda_1$ at
its upper bound. Next, this process is repeated for $\lambda_2$, and we
also find that $\lambda_2$ can be fixed at its upper bound.

Next, we consider proving that $\tau_3 = \ldots = \tau_{27} = 0$.  First note that
$\tau_i$ can be either zero or positive.  If $\tau_i > 0$, then $\lambda_i = p_i/(1-\beta)$
is at its upper bound and $\tau_i = x^T y^i - m$ in (\ref{eq:VaR MPEC}).
However, since $1 - \lambda_1 - \lambda_2 < p_i/(1-\beta)$, it follows that the LP
corresponding to $\tau_i > 0$ must be inconsistent for all $i=3,\ldots, 27$.
This is represented in Figure~\ref{F-tree} by the grey nodes. In
practice, the preprocessor in AMPL detects that these LPs are inconsistent,
and no solves are necessary.  Hence, we conclude that $\tau_3 = \ldots = \tau_{27} = 0$.
Solving the LP relaxation (\ref{eq:LP relax}) with these $\tau_i$ fixed at zero and also
$s_1 = s_2 = 0$ (since $\lambda_1$ and $\lambda_2$ are at their upper bounds),
we obtain a lower bound of $4.2652$,
which means that our candidate solution is globally optimal.

All in all, we solved only six LPs to prove the global optimality of the
candidate solution. The empty nodes are never solved, while the grey
nodes can be eliminated by the preprocessor.

For this example, we have compared the two lower bounds (\ref{eq:LP relax})
and (\ref{eq:convex relax}) with $\beta = 0.8$ and $\beta = 0.9$.
With $\beta=0.8$, the LP relaxation gives a lower bound of $m=0.61$,
which is poorer than the convex hull relaxation, which gives $m=1.24$.
However, with $\beta=0.9$, the LP relaxation (\ref{eq:LP relax})
gives the tighter bound with $m=3.48$, while (\ref{eq:convex relax})
gives only $m=2.45$.  Hence, we use the LP relaxation
(\ref{eq:LP relax}) in the above report.  (For this example, both
lower bounds would actually generate identical trees except for the
value at the root.)  We have also used the same branch-and-cut
procedure to verify global optimality in the problem with $\beta = 0.8$.
Apart from the fact that the details are similar, with $\beta = 0.8$,
the vector $x^{\scCVaR}$ obtained by solving the CVaR minimization problem
produced a VaR that is already globally optimal, as verified by the
branch-and-cut procedure.  In other words, the upper bounding refinement
is not needed in the case where $\beta = 0.8$; for this reason, we omit the
details. \hfill $\Box$

\section{Approximation by Smoothing}

In two pioneering papers \cite{CMangasarian95,CMangasarian96}, Mangasarian
and his then-Ph.D.\ student Chen developed a class of smoothing methods
for solving complementarity problems.  The basis of their methods is a family
of smooth functions that approximate the plus function $t_+$.  A summary of these
functions can be found in \cite[Subsection 11.8.2]{FPang03}.  In what follows,
we show how a smoothing approach can be applied to the VaR minimization problem.

Let $\rho_{\varepsilon}$ be any nonnegative-valued, twice
continuously differentiable, strictly convex function defined on the real
line such that $|\rho_{\varepsilon}^{\, \prime}(t)| \leq 1$ and
$\rho_{\varepsilon}^{\, \prime \prime}(t) > 0$ for all $t$,
and for some constant $c > 0$,
\begin{equation} \label{eq:uniform approx}
| \, t_+ - \rho_{\varepsilon}(t) \, | \, \leq \, c \, \varepsilon,
\epc \forall \, t \, \in \, \Re,
\end{equation}
for all $\varepsilon > 0$ sufficiently small.  The latter approximating
property has several consequences; among these, we have
\begin{equation} \label{eq:asymptotic rho}
\lim_{t \to \infty} \, \rho_{\varepsilon}(t) \, = \, \infty
\epc \mbox{and} \epc
\lim_{t \to -\infty} \, \rho_{\varepsilon}(t) \, = \, 0.
\end{equation}
Two examples of such a smoothing function are
\[
\rho_{\varepsilon,1}(t) \, \equiv \, \varepsilon \, \log( 1 + e^{t/\varepsilon} )
\epc \mbox{and} \epc
\rho_{\varepsilon,2}(t) \, \equiv \, \frac{\sqrt{t^2 + 4 \varepsilon^2} + t}{2}.
\]
We leave it to the reader to verify that these functions satisfy the assumed
properties.

In general, given a smoothing function $\rho_{\varepsilon}$ with the cited
properties, consider the $\varepsilon$-approximate CVaR and VaR minimization problems:
\[
\left\{ \begin{array}{ll}
\mbox{minimize}   & \CVaR_{\beta,\varepsilon}(x) \\ [5pt]
\mbox{subject to} & x \, \in \, X
\end{array} \right\} \epc \mbox{and} \epc
\left\{ \begin{array}{ll}
\mbox{minimize}   & \VaR_{\beta,\varepsilon}(x) \\ [5pt]
\mbox{subject to} & x \, \in \, X
\end{array} \right\},
\]
where
\[
\CVaR_{\beta,\varepsilon}(x) \, \equiv \, \min_{m \in \Re} \, \left[ \, m + \frac{1}{1 - \beta} \,
\displaystyle{
\sum_{i=1}^k
} \, p_i \, \rho_{\varepsilon}(x^Ty^i - m) \, \right],
\]
and $\VaR_{\beta,\varepsilon}(x)$ is the unique minimizer in the definition
of $\CVaR_{\beta,\varepsilon}(x)$.  It follows from the assumed properties of $\rho_{\varepsilon}$
that the minimand in $\CVaR_{\beta,\varepsilon}(x)$ is a coercive function of $m$, for fixed $x$,
that is,
\[
\lim_{m \to \pm \infty} \, \left[ \, m + \frac{1}{1 - \beta} \,
\displaystyle{
\sum_{i=1}^k
} \, p_i \, \rho_{\varepsilon}(x^Ty^i - m) \, \right] \, = \, \infty.
\]
(Use the second limit in (\ref{eq:asymptotic rho}) to show the case $m \to \infty$ and
the first limit and the nonexpansiveness of $\rho_{\varepsilon}$, which in turn is implied
by the condition $|\rho_{\varepsilon}^{\, \prime}(t)| \leq 1$ for all $t$, to show the other case
$m \to -\infty$.)  Therefore, $\CVaR_{\beta,\varepsilon}(x)$
is a well-defined, finite scalar.  Moreover, by the strict convexity of $\rho_{\varepsilon}$,
it follows that $\VaR_{\beta,\varepsilon}(x)$ is uniquely defined.  The latter is a significant
difference from the original $\CVaR_{\beta}(x)$ where the minimizing set ${\cal M}_{\beta}(x)$
is often not a singleton.  Another important consequence with using the smooth function
$\rho_{\varepsilon}$ is that $\VaR_{\beta,\varepsilon}(x)$ can be characterized as the unique
scalar $m$ that satisfies the smooth equation:
\[
1 - \beta \, = \, \sum_{i=1}^k \, p_i \, \rho_{\varepsilon}^{\, \prime}(x^Ty^i - m).
\]
This implies, by the implicit-function theorem, that the $\VaR_{\beta,\varepsilon}(x)$ is a
continuously differentiable function of $x$ with gradient given by
\[
\nabla \VaR_{\beta,\varepsilon}(x) \, = \,
\frac{ \displaystyle{
\sum_{i=1}^k
} \, p_i \, \rho_{\varepsilon}^{\, \prime \prime}(x^Ty^i - \VaR_{\beta,\varepsilon}(x)) \, y^i}{
\displaystyle{
\sum_{i=1}^k
} \, p_i \, \rho_{\varepsilon}^{\, \prime \prime}(x^Ty^i - \VaR_{\beta,\varepsilon}(x))} .
\]
The upshot of these properties is that $\VaR_{\beta,\varepsilon}(x)$ has much nicer analytic
properties than $\VaR_{\beta}(x)$; furthermore, the $\varepsilon$-approximate VaR minimization
problem is a smooth, albeit still nonconvex, linearly constrained nonlinear program in the
sole variable $x$.  As such, there are a host of efficient algorithms that one can use for
computing the minimum value (to be precise, stationary values) of the
$\varepsilon$-approximate value-at-risk.

An important question that arises is what happens to the convergence of the
$\varepsilon$-approximation problems as $\varepsilon \downarrow 0$.  Although
such a question has been partially studied in a general context (see, e.g., \cite{Kall86}),
we give a self-contained treatment to such a convergence issue for our special problem.
For this purpose, we establish a preliminary boundedness lemma.

\begin{lemma}  \label{lm:convergence of m} \rm
Let $\{\varepsilon_{\nu}\}$ be a sequence of sufficiently small positive
scalars, and let $\{x^{\nu}\}$ be an arbitrary sequence of vectors in $X$,
both of which are necessarily bounded.  The sequence $\{m_{\nu}\}$, where
$m_{\nu} \equiv \VaR_{\beta,\varepsilon_{\nu}}(x^{\nu})$
for every $\nu$, is bounded.  Moreover, if the pair $(m_{\infty},x^{\infty})$ is the
limit of a convergent subsequence $\{(m_{\nu},x^{\nu}) : \nu \in \kappa\}$ corresponding
to a sequence $\{\varepsilon_{\nu}\}$ of positive scalars tending to zero, then $m_{\infty}$
is an element of ${\cal M}_{\beta}(x^{\infty})$.
\end{lemma}

\noindent {\bf Proof.}  We have, for any $m \in \Re$ and any $\nu$,
\[
m_{\nu} + \displaystyle{
\frac{1}{1-\beta}
} \, \displaystyle{
\sum_{j=1}^k
} \, p_j \, \rho_{\varepsilon_{\nu}}((x^{\nu})^Ty^j - m_{\nu}) \, \leq \,
m + \displaystyle{
\frac{1}{1-\beta}
} \, \displaystyle{
\sum_{j=1}^k
} \, p_j \, \rho_{\varepsilon_{\nu}}((x^{\nu})^Ty^j - m),
\]
which implies, by the uniform approximation property (\ref{eq:uniform approx}),
\[
m_{\nu} + \displaystyle{
\frac{1}{1-\beta}
} \, \displaystyle{
\sum_{j=1}^k
} \, p_j \, ((x^{\nu})^Ty^j - m_{\nu})_+ \, \leq \,
m + \displaystyle{
\frac{1}{1-\beta}
} \, \displaystyle{
\sum_{j=1}^k
} \, p_j \, ((x^{\nu})^Ty^j - m)_+ + \frac{2 \, c \, \varepsilon_{\nu}}{1 - \beta}.
\]
Since the right-hand side is bounded, the boundedness of the sequence $\{m_{\nu}\}$
follows readily.  Moreover, the second assertion of the lemma also follows easily
from the last inequality.  \hfill $\Box$

\gap

We introduce the notation for our next result, which addresses the main convergence
issue of the $\varepsilon$-smoothing procedure.  Specifically, let
$\{\varepsilon_{\nu}\}$ be an arbitrary sequence of positive scalars converging
to zero.  For each $\varepsilon > 0$, let $x^{\varepsilon}$ be a (globally) optimal
solution of the smooth optimization problem:
\begin{equation} \label{eq:varepsilon VaR}
\begin{array}{ll}
\mbox{minimize}   & \VaR_{\beta,\varepsilon}(x) \\ [5pt]
\mbox{subject to} & x \, \in \, X,
\end{array} \end{equation}
which is clearly equivalent to
\[ \begin{array}{ll}
\mbox{minimize}   & m \\ [5pt]
\mbox{subject to} & x \, \in \, X \\ [5pt]
\mbox{and} & 1 - \beta \, = \,
\displaystyle{
\sum_{i=1}^k
} \, p_i \, \rho_{\varepsilon}^{\, \prime}(x^Ty^i - m).
\end{array} \]
Let $m_{\varepsilon} \equiv \VaR_{\beta,\varepsilon}(x^{\varepsilon})$.

\begin{proposition} \label{pr:convergence of approx} \rm
Suppose that the VaR minimization problem has a minimizer $x^{\scVaR}$ such that
${\cal M}_{\beta}(x^{\scVaR})$ is the singleton $\{m_{\scVaR}\}$.  It holds that
\begin{equation} \label{eq:limit approx VaR}
\lim_{\varepsilon \downarrow 0} \, m_{\varepsilon} \, = \, \min_{x \in X} \, \VaR_{\beta}(x).
\end{equation}
Consequently, for any sequence of positive scalars $\{\varepsilon_{\nu}\}$
converging to zero, every accumulation point of $\{x^{\varepsilon_{\nu}}\}$ is a minimizer
of $\VaR_{\beta}(x)$ on $X$.
\end{proposition}

\noindent {\bf Proof.}  Let $\{\varepsilon_{\nu}\}$ be an arbitrary sequence of positive
scalars tending to zero.  We have, for every $\nu$,
\[
m_{\varepsilon_{\nu}} \, \leq \, \VaR_{\beta,\varepsilon_{\nu}}(x^{\scVaR}).
\]
By Lemma~\ref{lm:convergence of m}, every accumulation point of the
right-hand side is an element of ${\cal M}_{\beta}(x^{\scVaR})$,
which by assumption is the singleton $\{m_{\scVaR}\}$.  Consequently,
the right-hand side converges to $m_{\scVaR}$ as $\varepsilon_{\nu} \downarrow 0$.
If $(m_{\infty},x^{\infty})$ is the limit of a convergent subsequence
$\{(m_{\varepsilon_{\nu}},x^{\varepsilon_{\nu}}) : \nu \in \kappa\}$,
then $m_{\infty} \in {\cal M}_{\beta}(x^{\infty})$,
and, from the above inequality,
\[
m_{\infty} \, \leq \, m_{\scVaR} \, = \, \VaR_{\beta}(x^{\scVaR});
\]
this shows that $x^{\infty}$ is a minimizer of $\VaR_{\beta}(x)$ on $X$ and
also establishes the desired limit (\ref{eq:limit approx VaR}).  \hfill $\Box$

\gap

Proposition~\ref{pr:convergence of approx} is theoretically very desirable; its
practical drawback is that there is no guarantee that a globally optimal solution
to (\ref{eq:varepsilon VaR}) can be computed.
We have implemented the two smoothing functions $\rho_{\varepsilon,1}$ and
$\rho_{\varepsilon,2}$ using $\varepsilon=10^{-3}$ for the example of the previous
section.  To solve the smoothed problem, we used five state-of-the-art NLP solvers:
filter, knitro 3.0, loqo 6.06, minos 5.5, and snopt 6.6-1, which are all available
on the NEOS server \cite{CzyzMesnMore:98,NEOS}.
The MPEC solver required 5 iterations to produce an upper bound for the example
on hand, which turns out to be globally optimal.  In contrast,
the NLP solvers fail for $\rho_{\varepsilon,1}$
because the exponentials cannot be evaluated or blow up during the computation.
The situation is slightly better for the square root formulation.  Filter produces
an optimal solution $m=4.265827$ (which is slightly higher than the minimum VaR of
4.2652); knitro produces a local optimum $m=4.74277$.
All other solvers fail to produce a feasible point (minos and snopt), while
loqo fails because it reached its iteration limit.

\section{Conclusion}

In this paper, we have investigated the minimization problem of the VaR as a nonconvex
LPEC and developed bounding schemes that can be used to verify the global optimality
of a candidate feasible solution.  We have also established the convergence of a
smoothing approach to compute an approximate VaR.  Whereas the VaR minimization problem
is special (and yet important in its own right), we maintain that the
bounding schemes can be extended to more general LPECs, and possibly even to other
``convex'' MPECs, namely, MPECs whose only nonconvexity is the complementarity constraint.
Indeed, the extension of the upper bounding scheme is fairly straightforward; it is
the lower bounding scheme that is very much problem dependent.  Nevertheless, we believe
that for special classes of MPECs, tight lower bounds can be obtained, which can
then be used in a branch-and-cut scheme either for verifying the global optimality
of a candidate solution obtained from a local MPEC solver or for computing a globally optimal
solution to the problem directly.

A lesson we have learned from the computational experiments in this paper is that
while the NLP solvers are generally very robust, one still requires a proof such as
the one given in Subsection~\ref{subsec:branching optimality} to ascertain the
quality of the solutions they produce.  For MPECs, we believe that the time is now
ripe for combining existing local methods with some global branch-and-cut schemes
in order to obtain solutions with proven global optimality.

\gap
\noindent {\bf Acknowledgment.} Professor Stanislav Uryasev has kindly alerted us to his
paper with Larsen and Mausser \cite{LMUryasev02} in which CVaR-based algorithms are
proposed for solving the VaR-minimization problem.

\vfill
\begin{flushright}
\scriptsize
\framebox{\parbox{2.4in}{The submitted manuscript has been created
in part by the University of Chicago as Operator of Argonne
National Laboratory ("Argonne") under Contract No.\
W-31-109-ENG-38 with the U.S. Department of Energy.
The U.S. Government retains for itself, and others
acting on its behalf, a paid-up, nonexclusive, irrevocable
worldwide license in said article to reproduce,
prepare derivative works, distribute copies to the
public, and perform publicly and display publicly, by or on
behalf of the Government.}}
\normalsize
\end{flushright}

\end{document}